\documentclass [winedt,yap]{iitparc}
\usepackage{cite}
\usepackage{amsmath,amssymb,amsfonts,bm}
\usepackage{lscape}
\usepackage{floatfig,wrapfig,epsfig}
\usepackage{subfigure}
\usepackage{color}
\usepackage{psboxit}
\usepackage{rotating}
\usepackage{curves}
\usepackage{ulem}
\usepackage[mathscr]{eucal}
\RequirePackage{psfrag}
\RequirePackage{graphicx}

\begin{document}\normalem
%
\frontmatter          
%
%
\IssuePrice{25.00}%
\TransYearOfIssue{2006}%
\TransCopyrightYear{2006}%
\OrigYearOfIssue{2006}%
\OrigCopyrightYear{2006}%
\TransVolumeNo{67}%
\TransIssueNo{2}%
\OrigIssueNo{2}%
%
\mainmatter              
%
\setcounter{page}{311}
\CRubrika{CONTROL IN SOCIAL ECONOMIC SYSTEMS}
\Rubrika{CONTROL IN SOCIAL ECONOMIC SYSTEMS}
%
\newtheorem*{conj*}{Conjecture}{\bfseries}{\itshape}

\title{Strategies of Voting in Stochastic Environment:  \\ Egoism and
Collectivism\thanks{This work was supported in part by the Russian
Foundation for Basic Research, project no.~02-01-00614.}}

\titlerunning{Strategies of Voting}
\author{V.  I.  Borzenko, Z.  M.  Lezina, A.  K.  Loginov,\\
Ya. Yu. Tsodikova, and P.  Yu.  Chebotarev}

\authorrunning{Borzenko \lowercase{$et~al.$}}

\OrigCopyrightedAuthors{Borzenko, Lezina, Loginov, Tsodikova, Chebotarev}

\institute{Trapeznikov Institute of Control Sciences, Russian Academy of
Sciences, Moscow, Russia}
\received{Received February 1, 2005}
\OrigPages{154--173}

\maketitle 

\begin{abstract}
Consideration was given to a model of social dynamics controlled by
successive collective decisions based on the threshold majority procedures.
The current system state is characterized by the vector of participants'
capitals (utilities).  At each step, the voters can either retain their
{\it status quo\/} or accept the proposal which is a vector of the
algebraic increments in the capitals of the participants.  In this version
of the model, the vector is generated stochastically.  Comparative utility
of two social attitudes---egoism and collectivism---was analyzed.  It was
established that, except for some special cases, the collectivists have
advantages, which makes realizable the following scenario:  on the
conditions of protecting the corporate interests, a group is created which
is joined then by the egoists attracted by its achievements.  At that,
group egoism approaches altruism.  Additionally, one of the considered
variants of collectivism handicaps manipulation of voting by the
organizers.
\end{abstract}

\DOI{0093}

\section{introduction}

As was shown by A.V.~Malishevskii [1, pp.~92--95] at a turn of the 1970's,
for any distribution of riches in society one can formulate a number of
proposals such that each of them enriches 99\% members of the society, but
as a result of successive acceptance of these proposals each member of the
society becomes very poor.  Bluntly speaking, these proposals lie in
successive dispossession of all participants one after another.  At that,
part of the participant's property is divided between the rest of them
including the already dispossessed ones, another part is passed to the
person formulating the proposal.  This algorithm can be used in a number of
votings, provided that the participants vote according to their valuable
interests and completely ignore the interests of the others.  One can
readily see that the heart of the problem remains unchanged if instead of
defending their own interests the voters defend the interests of their
electors, clans, industries, countries (in the UN), and so on.  This is due
to the fact that if the voters are guided by any particular interests
dividing them into smaller groups, then voting becomes easily manipulatable
by those who have priority in formulating the proposals submitted to voting
(see also [2, 3]).

How to eliminate this manipulation?  It is difficult to restrict the powers
of the ``Presidium'' and the influence of the political backstage spin
doctors.  The passage from shorthanded protection by the voters of their
private interests to the universal principle of ``all for one and one for
all'' would be a reliable safeguard, but as a rule it is unlikely.  First,
the ``voters will not understand.'' Second, ``self likes itself best.''
From the practical point of view, the participant is interested only in
``all for one,'' and only if this ``one'' is himself.  Yet the participant
may hope never to be in the role of the ``one'' needing protection,
especially if he intends to remain truthful to the ``Presidium.'' Moreover,
he understands that if sometimes he will be in this position, the others
may betray him, no matter how arduously he defended their interests in the
past.  Therefore, the immediate utility of the everyday mission of ``being
for all'' for which he will have to pay a pretty penny is very doubtful for
him.

However, it is not always the case that the collective decisions are made
in the conditions of selfish or malevolent manipulation by the proposal
makers.  In essence, the agenda is often defined not by them, but rather by
external---natural, man-caused, economical, political, demographic, and so
on---phenomena.  In what concerns these phenomena, it is possible to
discuss only the general trends, rather than the times of occurrence and
the scope of the future events.  In the first approximation this reality
may be modeled by a sequence of stochastically generated proposals.  Each
proposal is advantageous to some participants and disadvantageous to some
other participants.  In the simplistic model, this proposal may be
identified with a vector of algebraic increments in some abstractly
understood capitals (utilities) of the participants.

A model of this sort is discussed below.  The main subject matter is
represented by the social attitudes of the voting participants.  The main
attitudes are egoistic and collectivistic.  Like in the Malishevskii model,
the egoistic participant supports any proposal providing an increment to
its capital.  There are also other participants making up (still one) group
and voting jointly for the proposals that are beneficial, at least
minimally, to the group as a whole.  They sacrifice to some extent their
private interests to the interests of their group.  This kind of
collectivism can be justly called the group egoism.  The present paper is
concerned mostly with the attitude---egoistic or collectivistic---that is
beneficial to an individual participant.  More specifically, consideration
is given to the dynamics of the mean capital of the egoists and the group
members.

The following hypothetical mechanism is of special interest.  Let under
certain conditions the collectivistic attitude be more advantageous.  Then
the egoists have an incentive to join the group.  If the group puts no
obstacles in this way and grows by including more and more members, then
its group egoism resembles more and more the actions in the interests of
all, that is, the basically altruistic interests.  It deserves noting that
in this decision model the entire society participates in decision making,
rather than the ``Parliament'' that represents it.

Under the assumption of utility of the collectivistic attitude, this
mechanism enables the following scenario based on the pragmatic interest
and not on the altruistic ideals of the participants.  Let some
participants make up a group taking upon itself the responsibility to vote
not for their own interests, but for those of their group.  Formulation of
the latter is, of course, a problem {\it per se\/}.  The group is open to
new members who join it little by little as they see that the capital of
its participants on the average grows faster than that of the egoists.  The
conditions under which this scenario is feasible are the same as those
under which the collectivist attitude is preferable to the egoist one.
These conditions which are defined by certain values of the model inputs
are the subject matter of our study.

Why we discuss only a single group whereas in practice we mostly observe
confrontation of more than one group?  This is due to the fact that in the
environment where collectivism brings the best results it is more
beneficial to the groups to unite, rather than to compete.  In the
politics, the main obstacle to this is represented by the ideological
differences, ambitions of the leaders, and the desire to find its own
political ``niche'' (see, for example, [4]).  Somehow or other, but the
mechanisms of interaction of several groups within the framework of a model
undoubtedly deserve investigation; and this investigation is projected.  We
also plan to consider socially oriented groups which support their poorest
members in order to prevent their ruin, the model variants where the
increments of capital depend on their current values, the variants using
mechanisms for collection of taxes and ``party dues'' within the groups,
and so on.  At the same time, we do not plan to model purely economic
mechanisms of reproduction and drain of capital because we aim at analyzing
the social, rather than economic phenomena.  We again repeat that in the
model the ``capital'' is meant in the most general (like ``utility'') and
not economic sense.

\section{model of social dynamics}

\subsection{Basic Model}

There are $n$ participants among which $n_e$ are egoists and $n_g=n-n_e$
belong to the group.  It will be convenient in what follows to define the
relation of the participants by the parameter $\beta=n_e/2n$, half of the
fraction of egoists among all participants.  The ``society'' coincides with
the set of participants.  Its state at each instant is described by the
$n$-dimensional vector of capitals whose $i$th component is a real number
characterizing the capital of the $i$th participant.  Defined is the vector
of initial capitals; by the ``proposal of the environment'' is meant the
vector $(d_1,\dots,d_n)$ of capital increments, where $d_i$ is the
algebraic increment in the capital of the $i$th participant according to
the proposal of the environment.

It is assumed below as in [5--7] that the proposals of the environment are
generated stochastically and their distribution remains unchanged from step
to step.  Namely, we assume that $d_i$ is a normally distributed random
variable whose mean value and variance are denoted, respectively, by $\mu$
and $\sigma^2$.  The values $d_1,\dots,d_n$ are assumed to be independent
in the aggregate.

We denote by $\xi_e$ and $\xi_g$ the random variables representing,
respectively, the fractions of those egoists among all egoists and those
group members among all group members to whom a random proposal of the
environment provides a positive increment of the capital.  These random
variables take values over the interval $[0,1]$.  We do not rule out the
possibility of zero increment and at calculating $\xi_e$ and $\xi_g$ use
the coefficient $0{.}5$ to take into account the numbers of participants
getting the zero increment.

At each step the participants learn the next proposal and vote in line with
their principles.  If the $i$th participant is egoist, he votes for a
proposal if and only if $d_i>0$; if $d_i<0$, he votes against; and if
$d_i=0$, he abstains from voting (gives half-vote ``for'' and half-vote
``against'').  All group members vote identically in compliance with the
group principle of voting.  In the present paper, we consider two
principles, ``A'' and~``B''.

{\it The group votes ``for'' a proposal $(d_1,\dots,d_n)$ if and only
if\dots}

--~{\bf Principle A}:  \dots {\it as the result of accepting it the number
of group members getting a positive increment in the capital exceeds the
number of group members getting the negative increment\/}:  $\xi_g>0{.}5$;

--~{\bf Principle B}:  \dots {\it the sum of the increments in the capitals
of group members is positive\/$:$\ $\sum d_i>0$ $($the sum is taken over
the participants of the group$)$\/}.

We denote by $\xi$ the fraction of participants {\it voting\/} for the
given random proposal.  Since the group votes jointly,
\begin{gather}\label{eta}
\xi  =  \left\{%
\begin{array}{ll}
    (n_e\xi_e+n_g)/n & \text{if the group supports the proposal,} \\
    n_e\xi_e/n,        & \text{otherwise,} \\
\end{array}%
\right.
\\\nonumber
     =  \left\{%
\begin{array}{ll}
    2\beta\xi_e+(1-2\beta) & \text{if the group supports the proposal,} \\
    2\beta\xi_e,           & \text{otherwise.} \\
\end{array}%
\right.
\end{gather}

Decision making follows the ``$\alpha$-majority'' procedure according to
which a proposal is accepted if and only if $\xi>\alpha$, where $\alpha
\in[0,1]$ is the {\it decision threshold}.  Consideration will be given not
only to the voting thresholds $\alpha \geqslant 0{.}5$, but also to $\alpha
<0{.}5$ which are used in the practice of voting for various initiatives
such as organization of a new parliamentary group, sending a letter of
inquiry to the Constitutional Court, initiating a referendum, putting a
question on the agenda, and so on.  Such initiative decisions often play in
the social life a role not smaller than the majority decisions.

If a proposal is accepted, then the corresponding vector of increments
$(d_1,\dots,d_n)$ is added to the current vector of capitals; otherwise,
the latter remains unchanged.  Passage is made to the next step where a new
proposal is considered.  Thus, the {\it voting trajectory} is constructed.
We are going to consider the dynamics of the mean capital of the egoists
and the group on such trajectories and again emphasize that in the model at
hand the environmental proposals are random, that is, consideration is
given to the utility of the egoistic and collectivistic attitudes under
stochastic uncertainty, rather than to deliberate manipulation of voting by
its organizers.  Nevertheless, the problem of manipulatability will not
drop out of sight.  So, one can note straightway that in fact the voting
principle~A does not contribute to solving this problem.  Indeed, if the
group demonstrates its efficiency and all egoists join it, then the
situation will return to the initial one where no group existed at all:
the group will continue to vote jointly exactly for those decisions for
which the simple majority would vote if they were egoists.  It is namely
this voting that is most readily manipulatable.  Nevertheless, it is of
interest to compare the dynamics of the mean capitals of the egoists and
the group guided by principle~A under stochastic uncertainty.  Principle~B
resembles in many respects principle~A, but excludes---in the case where
the group coincides with the entire society---manipulation by the
organizers as described by A.V.~Malishevskii.  Indeed, according to
principle~B the group supports only those proposals that replenish its
common stock.  Therefore, after a series of such decisions the capitals of
all participants cannot decrease.

\subsection{Additional Options of the Model}

It is of interest to consider dynamics of participant's ruin.  To analyze
this phenomenon, provided was a variant of the model where the participant
with the capital falling down to a negative value ``leaves the field,'' is
disregarded in the new proposal of the environment, and does not vote.
Additionally, the model provides as option the possibility for egoists to
join or leave the group.  Two types of conditions for joining or leaving
the group were considered.  (i) The egoist is ready to join a group if its
capital remains smaller than the mean capital of the group members during
$s_1$ successive steps, where $s_1$ is a parameter.  Correspondingly, a
member of the group is ready to leave it if its capital is smaller than the
mean capital of the egoists over $s_2$ moves; $s_1$ and $s_2$ can differ.
(ii) The same comparison is carried out for the capital {\it increments\/}
rather than for the capitals themselves.  For the transitions to be
smoother and less ``mechanistic,'' it is assumed that if the condition for
the first or second type of transition is satisfied, then the transition
does not occur of necessity but with a certain probability which is the
model parameter.

\section{some examples}

The present paper considers the basic model which disregards ruin, joining
the group, and exit from it.  Some regularities of these phenomena were
described in brief in [6].  In the following examples the number of
participants is 200, the rms deviation of the capital increment is
$\sigma=10$, and the group adheres to principle~B.  We consider first the
case of the neutral environment ($\mu=0$), the group of 50\% of all
participants ($2\beta=0{.}5$), and the decision threshold 50\% \ ($\alpha
=0{.}5$).  Let the initial capital of all participants be $a=700$.  The
typical dependence of the capitals of the group members and the egoists
vs.~the step number%
\footnote[2]{One should not assume that hundreds of steps are always
required for appreciable changes in the capitals.  The parameters are
deliberately taken here in such a way that the changes are slow and the
height of the graph steps is insignificant on purpose not to dim the
general trend by random fluctuations.}
on a voting trajectory is depicted in Fig.~1a.

The mean capital of the egoists is practically time-independent, and that
of the group grows uniformly.  Variation of the voting threshold within
rather wide limits does not affect the picture.  Is it possible that
simultaneously the capitals of one category of the participants grow and
those of the other, diminish?  This example is shown in Fig.~1b where the
environment is unfavorable ($\mu=-1$), the group is small
($2\beta=0{.}92$), and the voting threshold is much higher than one half
($\alpha =0{.}48$).  Does the fact that decisions may be made against the
opinion of majority play here the key role?  No, an increase in the
threshold up to $\alpha =0{.}5$ entails no essential changes (Fig.~1c).
One more example of similar dynamics arises under an unfavorable
environment ($\mu=-1$), overwhelming majority of the group
($2\beta=0{.}08$), and extremely low decision threshold $\alpha =0{.}07$
(Fig.~2a).

\begin{figure}[t] 
\centering{\includegraphics[clip]{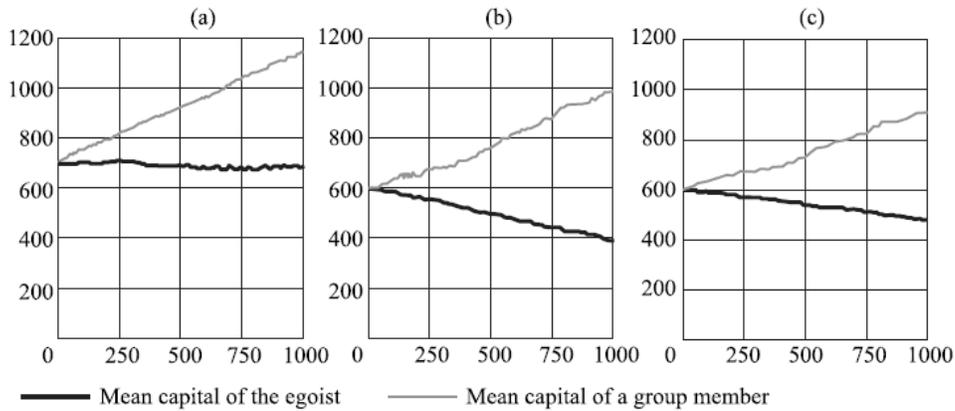}}
\vspace{-1.0em}\caption{The mean capital of egoists and group members vs.  the number of the step.  Two hundred
participants, $\sigma=10$, principle~B.  {\it a\/}) $\mu=0 $; $2\beta=0{.}5 $; $\alpha =0{.}5 $; {\it b\/})
$\mu=-1 $; $2\beta=0{.}92$; $\alpha =0{.}48$; {\it c\/}) $\mu=-1 $; $2\beta=0{.}92$; $\alpha =0{.}5 $.}\label{fig1}
\end{figure}

\begin{figure}[t] 
\centering{\includegraphics[clip]{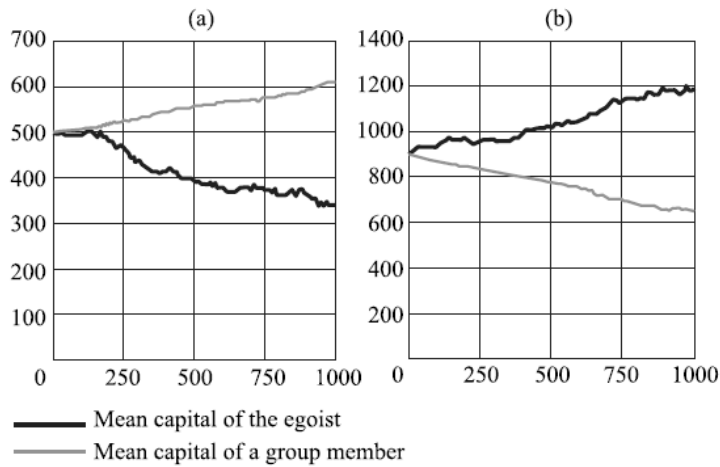}}
\vspace{-1.0em}\caption{The mean capital of egoists and group members vs.  the number of the step.  Two hundred
participants, $\sigma=10$, principle~B.  {\it a\/}) $\mu=-1 $; $2\beta=0{.}08$; $\alpha =0{.}07$; {\it b\/})
$\mu=-1 $; $2\beta=0{.}08$; $\alpha =0{.}04$.}\label{fig2}
\end{figure}

\begin{figure}[t] 
\centering{\includegraphics[clip]{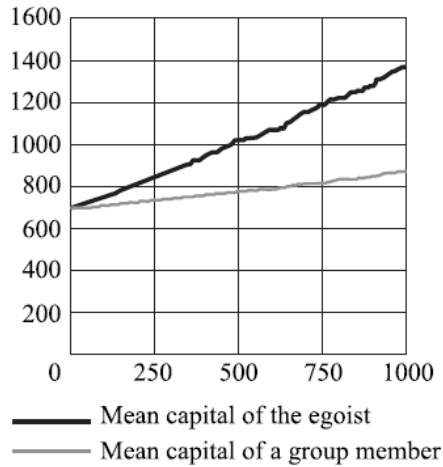}}
\vspace{-1.0em}\caption{The mean capital of egoists and group members vs.  the number of the step.  Two hundred
participants, $\sigma=10$, principle~B. $\mu=0{.}5$; $2\beta=0{.}08$; $\alpha =0{.}97$.}\label{fig3}
\end{figure}


At first sight it is surprising that even a greater reduction in the
decision threshold ($\alpha =0{.}04$) makes the curves practically change
places, the capital of the egoists grow and that of the group decrease
(Fig.~2b).  Finally, we discuss an example of high decision threshold.  For
$\alpha =0{.}97$, very large group ($2\beta=0{.}08$), and the favorable
environment ($\mu=0{.}5$), the egoists also are ahead (Fig.~3):  the mean
capital of the group member grows much slower than that of the egoists.
Consideration of these few examples convinces us that the regularities of
the social dynamics corresponding to the given model are not obvious;
therefore, it is of interest to analyze them by means of mathematical
tools.

\section{on random variables in the model}

In the given model, the proposal of the environment is the vector
$(d_1,\dots,d_n)$ of the capital increments of all participants, where
$d_1,\dots,d_n$ are independent random variables with the distribution
$N(\mu,\sigma^2)$.  The corresponding one-dimensional density and the
distribution function will be denoted by $f_{\mu,\sigma}(\cdot)$ and
$F_{\mu,\sigma}(\cdot)$;\ $f(\cdot)$ and $F(\cdot)$ denote the density and
the distribution function of the normal distribution with the center $0$
and variance $1$;\ $M(\eta)$ and $\sigma(\eta)$ denote the expectation and
the rms deviation of any random variable $\eta$ under consideration.

Each egoistic participant $i$ votes for a proposal if and only if $d_i>0$.
The probability of this event is as follows:
\begin{gather}\label{p}
p=P\{d_i>0\}=1-F_{\mu,\sigma}(0)=F\left(\frac{\mu}{\sigma}\right);
\end{gather}
the probability of voting ``against'' is as follows:
\begin{gather}\label{q}
q=1-p=P\{d_i<0\}=F_{\mu,\sigma}(0)=1-F\left(\frac{\mu}{\sigma}\right)=F\left(-\frac{\mu}{\sigma}\right).
\end{gather}

According to the model, the probability that the participant abstains from
voting is zero because the normal distribution is continuous.  Therefore,
the voting of each egoistic participant is the Bernoulli test with the
parameter $p$.  Then, since the values $d_i$ are independent, the number of
egoists voting ``for'' is distributed {\it binomially} with the parameters
$n_e$ and $p$.  The mean value and the variance of this distribution are,
respectively, $pn_e$ and $pqn_e$.  Normalization by dividing by $n_e$
provides the aforementioned random variable $\xi_e$, the {\it fraction of
egoists voting ``for''}.  It has the mean $p$ and the variance
$n_e^{-2}n_epq=n_e^{-1}pq$.  The fraction $\xi_g$ of the group members
getting a positive increment in capital according to the random proposal of
the environment has the same distribution with the mean $p$ and the
variance $n_g^{-1}pq$,\ $\xi_e$ and $\xi_g$ being independent.

The confidence interval---usually symmetrical or centered at
$M(\eta)$---which the random variable $\eta$ hits with the probability
$0{.}995$ will be called for brevity the {\it concentration zone} of
$\eta$.  In the case of normal random variable with the parameters $\mu$
and $\sigma^2$, this interval lies inside the segment
$[\mu-3\sigma,\,\mu+3\sigma]$ where about $99{.}73\%$ of the normal
distribution are concentrated.  Exit of the random variable from the
concentration zone will be regarded as a highly improbable event.  In a
series of, say, $1000$ steps, such events happen several times, but make no
appreciable contribution to the mean indices, which justifies their
screening-out.  The following sections will be devoted to analysis of the
social dynamics under various model parameters.

\section{case of neutral environment}

We first assume that $\mu=0$, that is, the environment is neutral and the
distribution of its proposals is symmetrical about $0$.  At that,
$p=q=0{.}5$; the distributions of the random variables $\xi_e$ and $\xi_g$
are symmetrical about $0{.}5$.  We assume for the time being, unless the
contrary is allowed, that the group is guided by the voting principle~A.
The effects of ruin and passage of the participant from one category to
another are disregarded.

\subsection{Decision Threshold $\alpha =1-\beta$} \label{se_1-bb}

Let us assume that the decision threshold is set as $\alpha =1-\beta$ and
that $2\beta<2/3$.  Then the voices even of all egoists will be
insufficient to make decision.  It is necessary and sufficient that the
proposal be supported by the group and at leat one half of the egoists,
that is, that the events $\xi_e>0{.}5$ and $\xi_g>0{.}5$, be realized (by
default the group uses the voting principle~A).  As was noted above, the
probabilities of each of these events are $0{.}5$, and the events are
independent.  Therefore, their joint probability is $0{.}25$, and,
therefore, the asymptotic (for a great number of steps) value of the
fraction of accepted proposals will be $0{.}25$.

Let us consider now the dynamics of the mean capital of egoists and the
group members.  As it was just established, the proposal is accepted if and
only if it provides for a positive increment in capital for more than half
of the egoists and more than half of the group members.  The mean {\it
positive} increment coincides with the {\it mean magnitude
$\sigma\sqrt{2/\pi}$ of the deviation from the mean} for the normal
distribution with the parameters $\mu$ and $\sigma$.  Let us estimate the
{\it number} of these positive increments.  The variance of the binomially
distributed number of egoists who voted ``for'' is $pqn_e=0{.}25n_e$, and
the rms deviation is $0{.}5\sqrt{n_e}$.  As in the case of the positive
increment, it is desired to determine the {\it magnitude of mean
deviation}.  Since the symmetrical binomial distribution is well
approximated by the normal distribution even for a relatively small number
of tests (usually this approximation is used for a number of Bernoulli
tests exceeding ${9}({pq})^{-1}$), in order to pass from the rms deviation
to the mean deviation magnitude we make use of the same coefficient
$\sqrt{2/\pi}$ as for the normal distribution, that is, \ estimate the mean
deviation magnitude by $\sqrt{{0{.}5n_e}/{\pi}}$.  Then, under the
condition that more participants voted ``for,'' the mean deviation of the
number of ``fors'' over the number of ``againsts'' is estimated by
$\sqrt{{2n_e}/{\pi}}$.  By multiplying it by the estimated positive
increment (this value and the number of positive increments are
independent), we get that the total increment for the egoists is
$\dfrac{2\sigma}{\pi}\sqrt{n_e}$.  Then, $\dfrac{2\sigma}{\pi\sqrt{n_e}}$
is the estimate of the capital increment of one egoist.  Since on the whole
a quarter of proposals is accepted, the mean increment for each egoist
participant in one step is estimated as $\dfrac{\sigma}{2\pi\sqrt{n_e}}$.
If the total number of steps is $s$, then after this series $\dfrac{\sigma
s}{2\pi\sqrt{n_e}}$ is the estimate of the expected capital increment of
the egoist.  For example, if $\sigma=10$ and $n_e=50$, then the estimated
increment for one participant in $100$ steps is about $22{.}5$.  Similarly,
the estimated mean capital increment in a series of $s$ moves for a member
of the group adhering to principle~A is $\dfrac{\sigma s}{2\pi\sqrt{n_g}}$.
Therefore, the growth rate of the mean capital of each category is in
inverse proportion to the root of its quantity (according to the well-known
pseudo-scientific aphorism, the speed of a human group is also in inverse
proportion to its size).  In particular, for $n_g=n_e$ the capital dynamics
of the group and egoists is the same, and for $n_g\ne n_e$ the smaller
category is in a better position.  However, one must bear in mind that the
considered decision threshold $\alpha =1-\beta$ itself depends on the
relation between the quantities.

Let now the group adhere to principle~B and give all its votes for the
proposal of the environment if and only if the total increment in the
capitals of its members is positive.  For a neutral environment, this
condition, as that of principle~A, is satisfied with the probability
$0{.}5$.  Since the increments in the capitals of the group and the egoists
are independent, a quarter of the proposals is accepted asymptotically as
before.  How the group dynamics changes at that?  For an arbitrary proposal
of the environment, the mean increment in the capital of the group member
is normally distributed with the mean $0$ and variance $\sigma^2/n_g$.  The
mean increment in the capital of the group member, provided that it is
positive as required by principle~B, is equal to the {\it mean magnitude}
of deviation, that is, differs from the rms deviation by the coefficient
$\sqrt{2/\pi}$ and is equal to $\sigma\sqrt{\dfrac{2}{\pi n_g}}$.  By
estimating again the fraction of the accepted proposals as one quarter, we
obtain for a series of $s$ moves the mean increment of the group member
equal to $\sigma s\left/\sqrt{8\pi n_g}\right.$.  The ratio of this value
to that obtained for the group using principle~A is
$\sqrt{\pi/2}\approx1{.}25$.  This 25\% gain can be attributed to the fact
that in both cases the group supports half of the proposals, but among the
supported proposals there are such that, although providing a positive
increment to more than half of its members, they nevertheless provide a
negative mean increment.  Therefore, in the case of principle~A the mean
capital grows slower than in the case of principle~B guaranteeing support
of the group even to the proposals which provide positive total increment
for the group and increase the capital only of the minority.

\begin{note}
The assumption of $2\beta<\frac{2}{3}$ made at the beginning of this
section can be relaxed using the notion of concentration zone of the random
variable $\xi_e$.  Indeed, to disable the egoists to make a decision by
their own forces without approval by the group, it is sufficient to the
accuracy of a highly improbable event that the threshold $\alpha =1-\beta$
be above the right boundary of the concentration zone of $2\beta\xi_e$.
This boundary is estimated by $2\beta\left(M(\xi_e)+3\sigma
(\xi_e)\right)$, where $M(\xi_e)=0{.}5$, $\sigma\left(\xi_e\right)=0{.}5/
\sqrt{n_e}$.  By solving the inequality $2\beta\left(0{.}5+3\cdot0{.}5/
\sqrt{n_e}\right)< \alpha =1-\beta$, we get
\begin{gather}\label{Wea_co1}
2\beta <\frac{1}{1+3\sigma \left(\xi_e\right)}=\frac{\sqrt{n_e}}{\sqrt{n_e}
+1{.}5}.
\end{gather}
For example, for $n_e=225$ this condition provides $2\beta<10/11$, thus
relaxing the initial constraint $2\beta<2/3$.
\end{note}

\subsection{Decision Threshold $\alpha =\beta$} \label{se_bb}

Now we consider the ``mirror'' case of voting with the threshold $\alpha
=\beta$ and assume for a start that $2\beta<2/3$.  At that, sufficient is
not only approval of the proposal by the majority of egoists, but also by
the group.  If in the above case a conjunction of both conditions (support
by the group and the majority of egoists) was required, here it suffices to
satisfy their disjunction.  Since the environment is neutral, the
probability of none of the events is 1/4; consequently, the disjunction of
the conditions is satisfied in 3/4 of cases.  These 75\% are divided into
50\% where the majority of the egoists are ``for'' and 25\% where the
majority are ``against.'' These 25\% are fully symmetrical to half of the
first of 50\% and ``balance'' them in the sense of the mean capital
increment (the sum of two means is zero).  We apply to the remaining half
of 50\% the same reasoning as in the case of the threshold $\alpha
=1-\beta$ and obtain the same result:  the mean capital increment of the
egoist is $\dfrac{\sigma s}{2\pi\sqrt{n_e}}$ over a series of $s$ steps.
The same results are obtained for the group:  the mean capital increment is
$\dfrac{\sigma s}{2\pi\sqrt{n_g}}$ or $\dfrac{\sigma s}{\sqrt{8\pi n_g}}$ depending
on which principle, A or B, is used.  The only difference of the general
dynamics lies in the fact that the presence of ``counterbalancing''
decisions, that is, providing the mean total zero of the ``quarters''
(fractions of 25\%), increases the spread as compared to the case of the
threshold $\alpha =1-\beta$ where they are absent.  Now, both in the group
and among the egoists a stronger stratification in incomes is observed.
Therefore, the voting thresholds $\alpha =1-\beta$ and $\alpha =\beta$ that
are symmetrical about $0{.}5$ provide different (by the factor of three)
numbers of the accepted proposals and different spreads of the incomes in
the group and among the egoists, but the same dynamics of the mean values.
As will be shown below, in the case of neutral environment the situation
with the mean values is always the same for $\alpha =\alpha '$ and $\alpha
=\alpha ''$ if $\alpha '+\alpha ''=1$.

\begin{note}
As in the above section, we relax the initial constraint $2\beta<2/3$.  Now
we need that the approval by the group should suffice (to within a highly
improbable event) for accepting the proposal of the environment.  By
assuming that support of the proposal by the egoists does not diminish the
left boundary of the concentration zone of $\xi_e$ which is equal to
$M(\xi_e)-3\sigma \left(\xi_e\right)$, we obtain that the minimal total
support of the proposal as approved by the group is $2\beta\left(M(\xi_e)-
3\sigma \left(\xi_e\right)\right)+(1-2\beta)$ because $1-2\beta$ is the
fraction of the group among the participants.  Therefore, the condition for
sufficiency of support by the group for approval of the proposal looks like
$2\beta\left(M(\xi_e)-3\sigma\left(\xi_e\right) \right)+ (1-2\beta)>\alpha
=\beta$.  By substituting $M(\xi_e)=\beta$ and $\sigma\left(\xi_e\right)=
0{.}5/\sqrt{n_e}$, we get the same condition (\ref{Wea_co1}) as in the
preceding case.
\end{note}

\subsection{Other Values of the Decision Threshold}

Since the symmetrical binomial distribution is well approximated by the
normal distribution even for a relatively small number of tests, we still
estimate the concentration zone of the symmetrical random variable $\xi_e$
by the interval with the boundaries $\left(M(\xi_e)\pm3\sigma\left(\xi_e
\right)\right)$, where in this case $M(\xi_e) =p=0{.}5$ and $\sigma\left(
\xi_e\right)=\sqrt{pq/n_e}=0{.}5/\sqrt{n_e}$, that is, the interval with
the boundaries $0{.}5\left(1\pm3/\sqrt{n_e}\right)$.  Then, according to
(\ref{eta}), it is unlikely that the random variable $\xi$ misses the union
of the interval with the boundaries
\begin{gather}\label{int1}
\beta\pm3\frac{\beta}{\sqrt{n_e}}\quad \text{(if the group does not support
the proposal)}
\end{gather}
and the interval with the boundaries
\begin{gather}\label{int2}
1-\beta\pm3\frac{\beta}{\sqrt{n_e}}\quad \text{(if the group supports the
proposal),}
\end{gather}
For example, if $n_e=n_g=225$, then $\beta=1/4$ and the values of $\xi$ are
concentrated in the domain $[0{.}2;0{.}3]\cup[0{.}7;0{.}8]$.

Let us consider the five cases of localization of the decision thresholds
$\alpha $.

{\bf Zone~1.} $\alpha <\beta-3\tfrac{\beta}{\sqrt{n_e}}$.  Here, almost
always $\xi>\alpha$, and practically all proposals of the environment are
accepted.  Since the environment is neutral, the expectation of the mean
increment in capital in the series of $s$ steps is extremely close to zero
both for the egoists and the group members, and the rms deviation of this
value is $\sqrt{s/n_e}\sigma$ for the egoists and $\sqrt{s/n_g}\sigma$ for
the group.
\smallskip

{\bf Zone~2.} $\alpha >1-\beta+3\tfrac{\beta}{\sqrt{n_e}}$.  Then, almost
always $\xi<\alpha$, and actually the proposals of the environment are
never accepted.  At that, the capital of the participants does not vary,
but the extremely rare accepted proposals increase the capitals of both the
group members and the egoists (the increase of the latter is greater
because they have a higher decision threshold than the group).
\smallskip

{\bf Zone~3.} $\beta+3\tfrac{\beta}{\sqrt{n_e}}<\alpha<1-\beta-3
\tfrac{\beta}{\sqrt{n_e}}$.  For these decision thresholds, $\xi>\alpha$ is
satisfied to an accuracy of low-probable events if and only if the group
supports the proposal.  If the group uses the voting principles~A or~B,
then this happens on the average in one half of cases, which is twice as
frequent as in the case discussed in Sec.~\ref{se_1-bb}
(page~\pageref{se_1-bb}), the rest being the same because the expected
value of the increment in the capital of the group member in the $s$-step
series is $\tfrac{\sigma s}{\pi\sqrt{n_g}}$ for principle~A and
$\tfrac{\sigma s}{\sqrt{2\pi n_g}}$ for principle~B.  The egoists do not
affect the voting; therefore, for them the proposals of the environment are
neutral on the average and the expected increment in their mean capital in
an $s$-step series is slightly greater than zero (nevertheless, slightly
greater owing to the rare improbable events) and has the rms deviation
$\sqrt{s/2n_e}\sigma$.  This zone of variations of the threshold $\alpha$
exists if the condition $\beta+3\tfrac{\beta}{\sqrt{n_e}}<1-\beta-3\tfrac{
\beta}{\sqrt{n_e}}$, which is equivalent to the condition $2\beta<\tfrac{
\sqrt{n_e}}{\sqrt{n_e}+3}$, is satisfied.  For example, for $n_e=225$ the
zone exists for $\beta<5/6$.
\smallskip

{\bf Zone~4.} $1-\beta-3\tfrac{\beta}{\sqrt{n_e}}\leqslant \alpha
\leqslant1-\beta+3\tfrac{\beta}{\sqrt{n_e}}$.  With an increase of the
decision threshold $\alpha$ from $1-\beta-3\tfrac{\beta}{\sqrt{n_e}}$ to
$1-\beta$ (zone ``4a''), the egoists exerts more and more influence on the
decisions and consistently reject more and more proposals that are
unadvantageous for them on the average.  At that, the mean increment in the
capital (in what follows, simply ``mean increment'') of an egoist in an
$s$-step series grows from zero (see case~3) to $\tfrac{\sigma
s}{2\pi\sqrt{n_e}}$ (Sec.~\ref{se_1-bb}, page~\pageref{se_1-bb}).  At the
same time, the increment of the group member in an $s$-step series is halved
from $\tfrac{\sigma s}{\pi\sqrt{n_g}}$ (zone~3) to $\tfrac{\sigma
s}{2\pi\sqrt{n_g}}$ (for principle~A) or from $\tfrac{\sigma s}{\sqrt{2\pi
n_g}}$ to $\tfrac{\sigma s}{2\sqrt{2\pi n_g}}$ (for principle~B).

\begin{figure}[t] 
\centering{\includegraphics[clip]{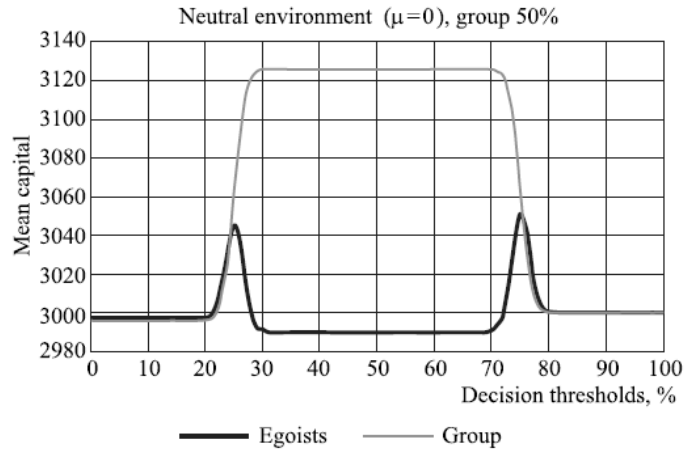}}
\vspace{-1.0em}\caption{The capitals averaged by the
categories and steps vs.  the threshold $\alpha $.  One realization, 450
participants, $\mu=0$, $\sigma=10$, principle~B, $2\beta=0{.}5$.}\label{fig4}
\end{figure}

With further growth of the threshold from $1-\beta$ to $1-\beta+3
\tfrac{\beta}{\sqrt{n_e}}$ (zone ``4b''), the mean increment of the group
member continues to vanish because the number of accepted proposals
vanishes.  For the egoists, the mean increment also diminishes because now
not all but only the most advantageous proposals are accepted.  An
interesting effect is observed here.  If in zone~4b $n_e=n_g$
($2\beta=0{.}5$) and the group adheres to principle~A, then the mean
increment is {\it higher\/} for the egoists rather than the group members
because the condition $\xi_g>0{.}5$ suffices to approve the proposal by the
group, whereas the fraction of the egoists' votes required for making a
decision is higher.  Therefore, the ``utility threshold'' of the accepted
proposals is higher for the egoists than for the group.  This effect was
mentioned when discussing zone~2, and one can readily see that it is
retained also when the group adheres to principle~B.
\smallskip

{\bf Zone~5.} $\beta-3\tfrac{\beta}{\sqrt{n_e}}\leqslant \alpha
\le\beta+3\tfrac{\beta}{\sqrt{n_e}}$.  Here, the dynamics of the mean
increments is a mirror reflection of that observed in zone~4.  When
$\alpha$ varies from $\beta-3\tfrac{\beta}{\sqrt{n_e}}$ to $\beta$ (zone
``5a''), less and less proposals for which $\xi\in(\beta-3
\tfrac{\beta}{\sqrt{n_e}},\,\beta)$ are accepted.  They are not supported
either by the group or by the majority of the egoists.  Therefore, in what
concerns the expected mean increments in the capital, they are unfavorable
to all, and a reduction in the number of such accepted proposals leads to
higher mean increments both for the group (principles~A and~B) and the
egoists.  Similar to the case of the ``mirror'' zone~4b, here also for
$n_g=n_e$ ($2\beta=0{.}5$) the mean capital increments are higher for the
egoists than for the group, which is accounted for by the fact that the
rejected proposals here are {\it just\/} unfavorable for the group
($\xi_g\leqslant0{.}5$) and {\it especially unfavorable} for the egoists
for which the stricter condition $\xi_e\leqslant \alpha /\beta<0{.}5$ is
satisfied.  Therefore, the egoists gain more from their rejection.

\begin{figure}[t] 
\centering{\includegraphics[clip]{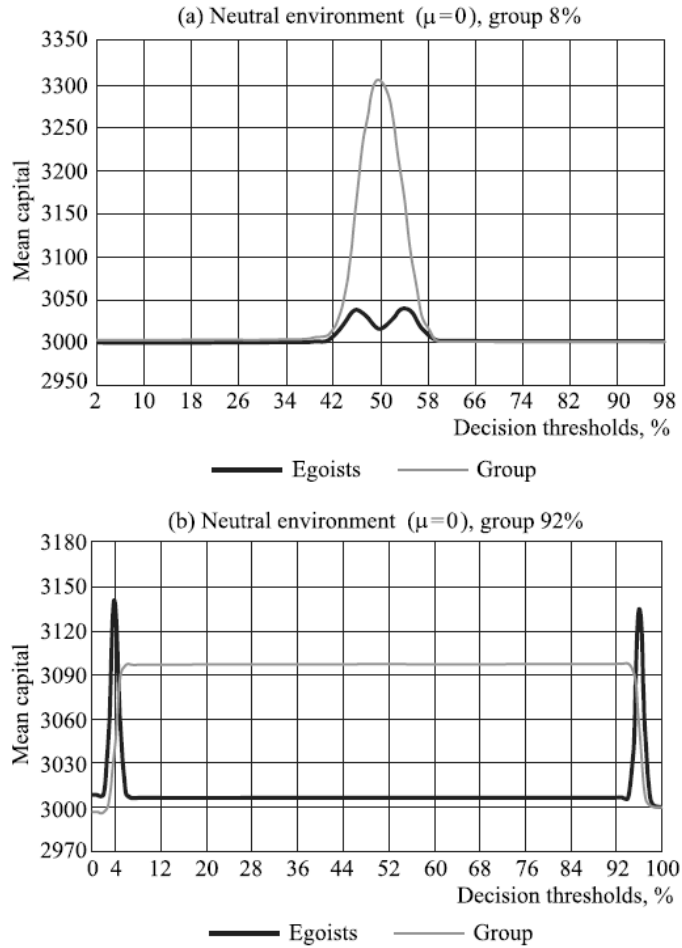}}
\vspace{-1.0em}\caption{The capitals averaged by the
categories and steps vs.  the threshold $\alpha $.  One realization, 450
participants, $\mu=0$, $\sigma=10$, principle~B.}\label{fig5}
\vspace{-0.5em}
\end{figure}


With further increase of $\alpha$ from $\beta$ to $\beta+3\tfrac{\beta}{
\sqrt{n_e}}$ (zone ``5b''), the set of the accepted proposals contracts
even more (their expected fraction decreases from 75\% to~50\%) at the
expense of the proposals rejected by the group but supported by the
majority of the egoists ($\xi>\beta$, consequently, $\xi_e>0{.}5$).
Therefore, the mean increment for the group continues to grow from the
values reached at $\alpha =\beta$ (Sec.~\ref{se_bb}, page~\pageref{se_bb})
to the value of zone~3, and the mean increment for the egoists decreases
from $\tfrac{\sigma s}{2\pi\sqrt{n_e}}$ (Sec.~\ref{se_bb},
page~\pageref{se_bb}) to zero.

The results of the computer-aided modeling are illustrated in Figs.~4 and~5
where the mean capitals of the egoists and the group members in a series of
$s=1000$ steps corresponding to various values of the decision threshold
are laid off on the vertical axis.  The initial capital of each participant
is 3000, the environment is neutral ($\mu=0$), $\sigma=10$, $n=450$; the
ruin and passage of the egoists to the group and exit from the group are
disregarded.  Averaging is carried out both over the participants and
steps.  By virtue of stationarity of the distribution of the capital
increments, the expectation of the {\it step mean} differs from the initial
value half as many as the expectation of the capital after the $s$th step.

These and the following figures depict the mean values for one realization
of the process which also reflect the spread in values about the
expectation, rather than the capital expectations (the data for it are
presented above, the general analytical expressions will be published
later).  For example, if the voting threshold lies in zone~3, the
expectation of the capital of egoists almost does not differ from their
initial capital, but their mean capital is appreciably (by 10 units)
smaller because of the spread in the graph (Fig.~4).

Figure~4 shows the case of as many egoists as there are members in the
group, $n_e=n_g$\ ($2\beta=0{.}5$).  The first graph of Fig.~5 refers to
the case of small group, $2\beta=0{.}92$, in which case the condition of
``Case~3'' cannot be met and no characteristic horizontal segment exists on
the curves of the mean capital of the egoists and the group.

In the example at hand, the value of $\beta$ is slightly higher than the
threshold of Notes~1 and~2 in which case for $\alpha =1-\beta$ a very small
fraction of decisions is made only by the egoists without approval of the
group and for $\alpha =\beta$ the proposal is not accepted despite support
of the group in a small fraction of cases.  The second graph of Fig.~5
shows the opposite case of a very large group, $2\beta=0{.}08$.  Here, the
greater part of the entire range of $\alpha$ lies in zone~3.  Since
$n_e<n_g$, for $\alpha =\beta\,$ and $\alpha =1-\beta\,$, the mean capitals
of the egoists are higher than those of the group members. Additionally,
owing to a low value of $n_e$, the capitals of the egoists have a
substantial spread for the values of $\alpha $ related to zones~1 and~3.
This spread accounts for the noticeable difference between the observed
mean capital of the egoists in zones~1 and~3 and the expected capital
coinciding with the initial capital (3000 units).

\subsection{Some More About the Zones where the Egoists Have Advantages}

In the case of neutral environment $n_e=n_g$, the above analysis shows that
in what concerns the expected increment in the capital the egoists have
advantage over the group members only in zones 4b and 5a.  In zone 4b, of
all proposals advantageous for the group the egoists take only those that
are most advantageous to them.  In the ``mirror'' case 5a, the egoists
block approval of the proposals that are {\it most\/} disadvantageous for
them, whereas the group cannot do so.

We note that for $n_e=n_g$ in these two zones the group can easily bring to
nothing the egoists' advantage by using the following procedure.  We define
for the group a special ``internal voting threshold'' $\alpha '$ depending
on $\alpha $ and $\beta$.
\begin{gather}\label{aa'}
\alpha '
 = \left\{%
\begin{array}{lll}
    \dfrac{1}{2}-\dfrac{\delta}{2\beta}, &\text{where}\quad\delta=\beta-\alpha  & \text{if}\quad\alpha
    <\beta;\\[3mm]
    \dfrac{1}{2}+\dfrac{\delta}{2\beta}, &\text{where}\quad\delta=\alpha -(1-\beta) & \text{if}\quad\alpha
    >1-\beta;\\[3mm]
    \dfrac{1}{2}  & & \text{if}\quad\beta\leqslant \alpha
    \leqslant1-\beta
\end{array}%
\right.=
\\\nonumber
 =
\left\{%
\begin{array}{ll}
    \dfrac{\alpha }{2\beta} & \text{if}\quad\alpha
    <\beta;\\[2mm]
    1-\dfrac{1-\alpha }{2\beta} & \text{if}\quad\alpha
    >1-\beta;\\[2mm]
    \dfrac{1}{2} & \text{if}\quad\beta\leqslant \alpha \leqslant1-\beta.
\end{array}%
\right.
\end{gather}

Let us consider principle~A$'$.

{\bf Principle A$'$}.  {\it The group votes for the proposal of the
environment if and only if in the case of its approval the fraction $\xi_g$
of its members getting a positive increment in capital exceeds the
threshold $\alpha '$ defined by $(\ref{aa'})$\/}.

Voting by principle A$'$ offers to the group the same possibilities of
influencing the decisions made in zones~4b and~5a as enjoyed by the
egoists.  Indeed, in zone 4b $\alpha =1-\beta+\delta$.  In the case of
accepting the proposal of the environment, $(1-2\beta)n$ votes are given by
the group; consequently, it is necessary that among all participants the
fraction of egoists supporting the proposal exceed $\alpha-(1-2\beta)=
\beta+\delta$.  Consequently, the condition $\xi_e>\frac{\beta+\delta}{2
\beta}=\frac{1}{2}+\frac{\delta}{2\beta}$ must be satisfied for approval of
the decision.  If the group establishes for itself the same threshold of
voting, that is, decision, it will be in the same position as the egoists:
they will have identical ``utility thresholds'' for the supported proposals
and, therefore, identical expected dynamics of the capital.

Similarly, in zone 5a $\alpha =\beta-\delta$.  To approve a proposal by the
efforts of egoists, that is, without participation of the group, it is
necessary that they provide more than $\alpha n=(\beta-\delta)n$ votes.
Therefore, among all egoists the fraction of those voting ``for'' must
exceed $\frac{\beta-\delta}{2\beta}=\frac{1}{2}-\frac{\delta}{2\beta}$.
Establishment of the same ``internal threshold'' smaller than $0{.}5$ puts
the group in the same conditions as the egoists:  if the votes collected in
the group exceed this threshold, then the proposal will be accepted for
sure, and this will happen as frequently as the votes of the egoists exceed the same threshold.  As the result, the expected dynamics of capital
in the group and among the egoists again will be the same.

If the number of egoists is smaller than one half, then under the above
changes in the intragroup threshold they will be again ahead in zones 5a
and 4b, but their advantage will decrease substantially.

Further reduction of the intragroup threshold in zone 5a and increase in
zone 4b%
\footnote[3]{This may be done in a most natural way by establishing an
intragroup threshold equal to the decision threshold $\alpha$.}
will offer advantage to the group over the egoists.  However, this
advantage will be relative, that is, the group really decreases its
increment in capital almost for all values of the threshold $\alpha $, but
the reduction of the egoists will be even greater.  If for the group to
``live better than the egoists'' is preferable just to just ``living
better,'' then it can reach this aim.

\section{cases of favorable and unfavorable environment}

The case of $\mu>0$ corresponds to the favorable environment, the case of
$\mu<0$, to the unfavorable environment.  Reasoning providing conclusions
about the nature of the dynamics of mean capitals are in this case similar
to those above.  Figure~6 illustrates the case of $n_e=n_g$\
($2\beta=0{.}5$), the rest of the parameters being the same as before.

\begin{figure}[t] 
\centering{\includegraphics[clip]{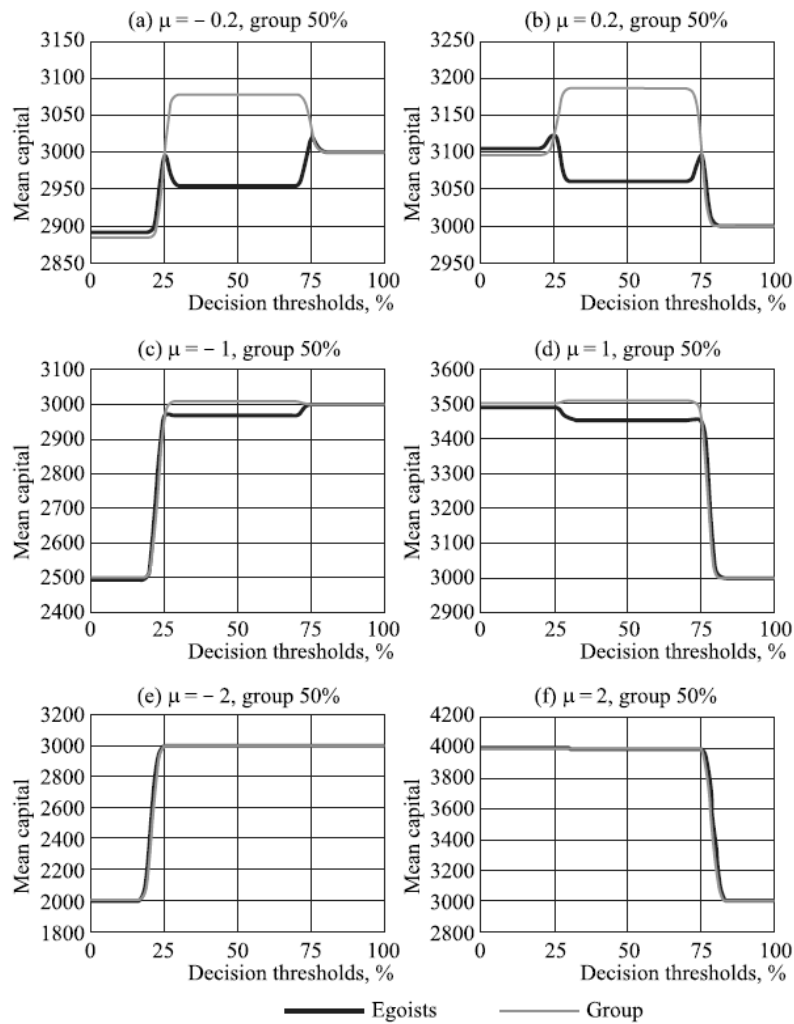}}
\vspace{-1.0em} \caption{Dependence of the
capitals averaged over the categories and steps on the threshold $\alpha$.  One realization, 450
participants, $\sigma=10$, $2\beta=0{.}5$, principle~B.}\label{fig6}
\end{figure}


For small deviations of $\mu$ from $0$ (Fig.~6a,b), the graphs have five
zones like those considered above.  In zone~2 (numeration as above), that
is, for the decision thresholds close to 1, the capital is equal to the
initial capital because the proposals are rejected.

On the contrary, in zone~1 almost all proposals are accepted.  Therefore,
the mean increment in capital after an $s$-step series is $\mu s$ (we recall
that the capitals of the participants are averaged on the graphs also over
the steps, therefore the values shown there are half as many) for the rms
deviation $\sqrt{s}\,\sigma$.

In zone 3, only those decisions are made that are supported by the group.
Let us estimate their frequency assuming that the group adheres to
principle~B.  The mean increment in capital of a group member in one step
$\widetilde{d}_{{\mathcal{G}}}$ has the distribution $N(\mu,(\sigma')^2)$,
where $\sigma'=\dfrac{\sigma}{\sqrt{n_g}}$.  The expectation of the
frequency of decisions made in zone 3 is equal to the probability of
positiveness of this value, that is,
\begin{gather}\label{P>0}
P\left\{\widetilde{d}_{{\mathcal{G}}}>0\right\}=F\left(\frac{\mu}{\sigma'}\right)=F(\mu'),
\end{gather}
where as before $F(\cdot)$ is the standard normal distribution function and
$\mu'$ stands for $\mu/\sigma'$.  Let us determine the expectation of
$\widetilde{d}_{{\mathcal{G}}}$, provided that it is positive.  After
integration we get
\begin{gather}\label{AveDgCo}
M \left(\widetilde{d}_{{\mathcal{G}}}\mid
\widetilde{d}_{{\mathcal{G}}}>0\right)
=\left(P\left\{\widetilde{d}_{{\mathcal{G}}}>0\right\}\right)^{-1}\int\limits_0^{\infty}
xf_{\mu,\sigma'}(x)\,dx =\mu+\sigma'\frac{f(\mu')}{F(\mu')}.
\end{gather}
Since the unconditional expectation of $\widetilde{d}_{{\mathcal{G}}}$ is
as follows:
\begin{gather}\label{AveDg}
M\!\left(\widetilde{d}_{{\mathcal{G}}}\right)
=P\left\{\widetilde{d}_{{\mathcal{G}}}>0\right\}M
\left(\widetilde{d}_{{\mathcal{G}}}\mid
\widetilde{d}_{{\mathcal{G}}}>0\right) =\sigma'f(\mu')+\mu
F(\mu'),
\end{gather}
the expectation of an increment in the capital of a group member after $s$
steps is as follows:
\begin{gather}\label{AveDgS} sM \left(\widetilde{d}_{{\mathcal{G}}}\right)
=s\left(\sigma'f(\mu')+\mu F(\mu')\right).
\end{gather}

The increment of the capital of an egoist in zone 3 after $s$ steps is
estimated by
\begin{gather}\label{AveDeS}
sM \left(\widetilde{d}_{{\mathcal{E}}} \right)
=sP\left\{\widetilde{d}_{{\mathcal{G}}}>0\right\}M
\left(\widetilde{d}_{{\mathcal{E}}}
\,|\,\widetilde{d}_{{\mathcal{G}}}>0\right) =s\mu F(\mu').
\end{gather}
Therefore, in zone 3 the estimated expectation of the difference in the
capitals of the group member and the egoist after $s$ steps is as follows:
\begin{gather}\label{DeltAve}
s\left(M \left(\widetilde{d}_{{\mathcal{G}}}\right)-M
\left(\widetilde{d}_{{\mathcal{E}}} \right)\right)
=s\sigma'f(\mu').
\end{gather}

As can be seen in Fig.~6, this value and $f(\mu')$ rapidly decrease with
increase in $|\mu|$:  as compared to the case of $|\mu|=0{.}2$, for
$|\mu|=1$ the graphs of the mean capital of the group members and the
egoists approach closely each other and actually fuse for $|\mu|=2$.  The
case of $n_e\ne n_g$ is illustrated in Fig.~7.

\begin{figure}[t] 
\centering{\includegraphics[clip]{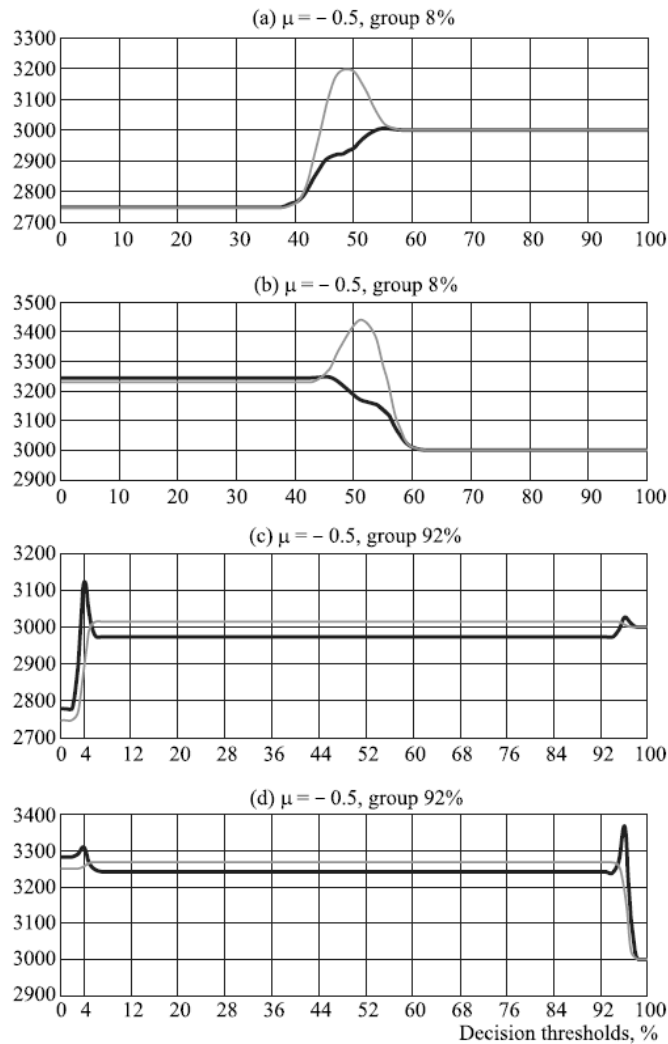}}
\vspace{-1.5em} \caption{Dependence of the capitals averaged in categories and steps vs. the threshold
$\alpha $.  One realization, 450 participants, $\sigma=10$, principle~B.}\label{fig7}
\end{figure}

The cases of $\beta=0{.}92$ and $\beta=0{.}08$ are shown here for
convenience of comparing with Fig.~5.  It goes without saying that the
regularities are the same:  with an increase in $\mu>0$, to the right of
zone 3 the graphs make a higher and higher ``step'' and a lower and lower
step to its left.  On the contrary, with a decrease of $\mu<0$, the graphs
have a higher ``step'' to the left of zone 3 and a lower one to the right
of it.  With increase in $|\mu|$, the graphs of the mean capitals of the
egoists and the group members draw together until almost complete fusion.

\section{summary}

Although the present paper models  the ``capital dynamics,''
the term
``capital'' should not be misleading.  It is not an economic model that is
considered above, but that concerning only the social and political
sciences.  Indeed, it does not cover reproduction of capital by investing
into production or using the financial tools, as well as changes in the
capital that are related with other factors independent of collective
decisions.  All this is unnecessary because it is required here to study
the dynamics defined by the mechanism of democratic voting with the
built-in basic social attitudes of egoism and collectivism, rather than the
economic or uncontrollable random mechanisms.  Therefore the term
``capital'' refers here to any countable resource (utility) controlled by
the collective decisions and having nothing to do with special regularities
of reproduction or waste.  In what follows, we summarize some conclusions
obtained by analyzing the model.

1.~The voting-dependent social dynamics is defined basically by the
decision threshold.  For higher thresholds (zone~2), the proposals are not
accepted and the {\it status quo\/} is retained.  For lower values
(zone~1), actually all proposals are accepted, the dynamics being identical
both for the group and the egoists and defined by the mean and the variance
of distribution of the proposals of the environment.  The zone boundaries
depend on the number of egoists and the group members, as well as on the
parameters of the environment.  With increase in the total number of the
participant and the fraction of egoists, zones~1 and~2 widen; for lower
values of these parameter they contract.  Inside the middle (zone~3), which
exists if the group is not too small, the dynamics is independent of the
decision threshold and the group has superiority over the egoists.  For the
neutral or a moderately unfavorable environment, it manages to realize
decisions that on the average are advantageous for its members.  At that,
the egoists do not influence the decisions, and for them the accepted
proposals do not differ from a random sample of the proposals of the
environment.  There are two zones on either side of zone~3 (zone~5 to the
left and zone~4 to the right) where the egoists influence the decisions.
There are two peaks (maxima) of the increment of the egoists' capital
which, however, vanish in the case of apparently favorable and apparently
unfavorable environments.  The group characteristics within these zones
vary monotonically from the values in the extreme zones~1 and~2 to the
value in the middle zone.  For the neutral environment, the expected
increment in the capital reached by the egoists at the two maxima---for
$\alpha =1-\beta$ and $\alpha =\beta$---are the same as for a group of the
same size for the same thresholds; for different sizes, the advantage
belongs to the smaller category.  If the egoists make 20\% of the total
number of participants, the increment in their capital reaches at the
maxima the value which the group has in the middle zone in the case of
voting by principle~A.  If the egoists make about 14\% of the participants,
then at the maxima they reach the value which the group has in the middle
zone in the case of the voting principle~B.  If the number of egoists is
still smaller, then for the neutral environment, the increment in their
capital at the maxima exceeds all increments reachable by the group.  For a
small group, zones~4 and~5 merge, and zone~3 degenerates.  At the
interfaces of zones~4 and~5 both the group and the egoists have maximal
increments in the capital, but the group maximum is higher.

2.~As for the scenario described in the Introduction---egoists join the
group and thus approach the group egoism to altruism,---it is absolutely
realistic.  The group is especially attractive for the mean values of the
decision threshold to which the simple-majority threshold $\alpha =0{.}5$
belongs.  For such thresholds, the group retains its attractiveness even if
it includes the majority of the participants; the egoists can have
advantage only for high and low thresholds.

3.~As was noted in Item~1 of this Summary, a smaller category of the
participants has advantage in zones~4 and~5.  Explanation of this curious
phenomenon is related with the law of large numbers:  since the variance of
the sample mean decreases with sample volume, the proposals that are ``very
good'' in the sense of the mean increment of capital are less frequent for
the larger category.  The ``very bad'' proposals are less frequent as well,
but this fact does not affect the dynamics because these proposals are
rejected by voting.

4.~Principle~B which better protects the group against the manipulations of
the organizers is preferable to principle~A also in the sense of the mean
increment in capital, which can be readily explained by the fact that it is
namely the positiveness of the mean (total) increment in capital, rather
than satisfaction of the majority of the group (as it is the case with
principle~A), that is declared by principle~B as the group utility.  For
example, in the case of the neutral environment, the ratio of the
increments in capitals of the groups adhering to principles~B and ~A is
$\sqrt{\pi/2}\approx1{.}25$.

5.~For changes in $\mu$ and $\sigma$ that retain their ratio $\mu/\sigma$,
the graphs of the capital increments extend/contract along the $Y$-axis in
proportion to $\sigma$.  Therefore, the impact of the parameter $\sigma$ on
the increments in capital is defined by $\mu$.  Namely, the passage from
$\sigma$ to $\sigma'=\rho\sigma$ ($\rho>0$) provides graphs extended by the
factor of $\rho$ that are obtained by passing from $\mu$ to $\mu/\rho$.

6.~The maxima of the increment in the egoists' capitals lie in the middles
of zones~4 and~5.  In the more distant, ``external'' parts of these zones,
the egoists have advantage over the group for smaller, equal, and even
somewhat higher their number.  The group can reduce or, for a sufficiently
high relative number of the egoists, even bring to nothing this effect by
passing to the voting principle~A$'$.  By changing its ``internal voting
threshold,'' the group can even achieve advantage over the egoists.  At the
same time, in the majority of cases it reduces its mean increment in
capital, but the mean increment in the egoists' capital reduces even more.
Is it advantageous to the group?  There is no unambiguous answer.  In the
social practice, the question ``What is more attractive, to live better
than before or to live better then the rest?'' always remains open.  It is
only clear that as a disciplined unit the group can choose an answer,
whereas the egoists do not have such a possibility.

\section{conclusions}

The paper analyzed the model of social dynamics defined by voting in the
stationary stochastic environment.  Time uniformity of the environment
parameter defines the specificity of the results and distinguishes the case
at hand from the situation of voting which the organizers try to
manipulate.  As was shown by analysis, for wider domains in the parameter
space, the group has a better dynamics of capital than the egoists, which
makes realistic the scenario where the egoists join the group and the group
egoism approaches altruism.  The narrow domains where the egoists have
advantage over the group were identified.  The main results of analysis
were interpreted in terms of the social sciences.

\revred{F.T. Aleskerov}

\begin{thebibliography}{7}
\itemsep0mm
\parsep0mm

\bibitem{Mirkin74} Mirkin, B.G., {\it Problema gruppovogo vybora} (Problem
of Group Choice), Moscow:  Nauka (Fizmatlit), 1974.

\bibitem{Aizerman81} Aizerman, M.A., Dynamic Aspects of the Voting Theory
(Review of the Problem), {\it Avtom.  Telemekh.}, 1981, no.~12,
pp.~103--118.

\bibitem{Chebotarev86} Chebotarev, P.Yu., Some Properties of Trajectories
in the Dynamic Problem of Voting, {\it Avtom.  Telemekh.}, 1986, no.~1,
pp.~133--138.

\bibitem{AleOrt} Aleskerov, F.T.  and Orteshuk, P., {\it Vybory.
Golosovanie.  Partii} (Elections.  Voting.  Parties), Moscow:  Akademiya,
1995.

\bibitem{Collective03} Borzenko, B.I., Lezina, Z.M., Lezina, I.B., {\it et
al.}, Model of Social Dynamics Defined by Collective Decisions and Random
Environmental Changes, in:  {\it II Mezhdunar.  konf.  po probl.  upr.\,
Tez.  dokl.} (II Int, Conf Control, Abstracts), Moscow:  IPU RAN, 2003,
p.~120.

\bibitem{Collective04a} Chebotarev, P.Yu., Borzenko, V.I., Lezina, Z.M.,
{\it et al.}, Model of Social Dynamics Controlled by Collective Decisions,
in:  {\it Tr.  In-ta probl.  upr.  RAN.} (Proc.  Int.  Control Probl.),
vol.~XXIII, Moskva, 2004, pp.~102--109.

\bibitem{Collective04b} Chebotarev, P.Ju., Borzenko, V.I., Lezina, Z.M.,
{\it et al.}, A model of Social Dynamics Governed by Collective Decisions,
in:  {\it Proc.  Int.  Conf.  ``Math.  Modelling of Social and Economic
Dynamics'' (MMSED-2004), June 23--25, 2004}, Moscow:  RSSU, 2004,
pp.~80--83.

\end{thebibliography}
\end{document}